
\input amstex
\magnification=\magstep1
\documentstyle{amsppt}
\topmatter
\title
            Spin polynomial invariants for Dolgachev surfaces
\endtitle
\author
             Stefan Bauer and Victor Pidstrigatch
\endauthor
\address
             Universit\"at Bielefeld, Germany and Insitute for Advanced Study,
Princeton, USA \hfill\break Steklov Institute, Moscow, Russia
\endaddress
\email
             bauer\@mathematik.uni-bielefeld.de  \,\,\,\,\,\,\,\,\,\,\,\,\,\,
             pidstrig\@alg.mian.su
\endemail
\thanks
The authors are grateful to the Mathematical Insitute in Oxford and to
the Sonderforschungsbereich 170 in G\"ottingen for their hospitality
during the research on this project. The first author was supported by a
Heisenberg stipend from the DFG
\endthanks
\endtopmatter

\define \lra{\longrightarrow}
\define \ra{\rightarrow}
\define \sA{{\Cal{A}}}\define \sC{{\Cal{C}}}
\define \sD{{\Cal{D}}}\define \sE{{\Cal{E}}}\define \sF{{\Cal{F}}}
\define \sG{{\Cal{G}}}\define \sH{{\Cal{H}}}
\define \sJ{{\Cal{J}}}\define \sL{{\Cal{L}}}
\define \sM{{\Cal{M}}}\define \sN{{\Cal{N}}}\define \sO{{\Cal{O}}}
\define \sR{{\Cal{R}}}
\define \sU{{\Cal{U}}}
\define \sV{{\Cal{V}}}

\define \fA{{\frak{A}}}\define \fB{{\frak{B}}}\define \fC{{\frak{C}}}
\define \fD{{\frak{D}}}\define \fF{{\frak{F}}}
\define \fG{{\frak{G}}}
\define \fK{{\frak K }}

\define \fT{{\frak{T}}}

\define \bC{{\Bbb C}}

\define \bP{{\Bbb P}}\define \bR{{\Bbb R}}

\define \bZ{{\Bbb Z}}
\define \coker{\operatorname{coker}}  \define \codim{\operatorname{codim}}
\define \can{\operatorname{can}}      
\define \Ext{\operatorname{Ext}}      
\define \Hilb{\operatorname{Hilb}}    \define \Hom{\operatorname{Hom}}
\define \im{\operatorname{im}}        \define \ind{\operatorname{ind}}
\define \length{\operatorname{length}}
     \define \Pic{\operatorname{Pic}}
\define \Proj{\operatorname{Proj}}    \define \rad{\operatorname{rad}}
       \define \Sign{\operatorname{Sign}}
    \define \Spin{\operatorname{Spin}}
\define \Spec{\operatorname{Spec}}    \define \supp{\operatorname{supp}}
\define \Sym{\operatorname{Sym}}      \define \w{\widetilde}
\document
\heading
I. Introduction
\endheading
The definition of Spin--polynomial invariants \cite{PT}, \cite{Ty} uses
the same machinery as Donaldson's invariants. One considers anti self dual
connections $a$ for which a coupled Dirac operator $\sD_a$ has nontrivial
kernel. For generic metrics the gauge equivalence classes of such
connections form a subspace of the moduli space of ASD--connections
which then can be used much
the same way as the moduli space itself to define invariants of the
differentiable structure of the underlying four dimensional manifold $X$.
To do all of this, one of course has to fix a $\Spin^\bC$--structure on $X$.

{}From the technical point of view these invariants appear to be better
accessible to algebraic geometric computations than the polynomial invariants
themselves.
This is because the corresponding bundles come with a section. In principle
this reduces the computation of the invariants to the study of the zeroes
of such a section and the associated extension groups.

In this paper we apply this seemingly simple recipe to a family of well known
algebraic surfaces, the Dolgachev surfaces, which have been studied
by a number of authors \cite{Do1}, \cite{FM1}, \cite{OV}, \cite{Ba1}. To
add some spice, there are a number of problems to tackle with. A first one
is the chamber structure: Since $b^2_+=1$, the invariants depend a priori
on the chosen metric and it seems hard to give a description of this
dependence in general.
However, for small Chern classes this dependence sometimes can be controlled.
In our case we get a diffeomorphism invariant, if we choose the Chern
classes $c_2=2$ and $c_1=K_S+2nk$. This invariant then is a polynomial
$$
q_S(n)=a(n)Q^2+b(n)Qk^2+c(n)k^4
$$
in the intersection form $Q$ and a primitive class $k$, which is a rational
multiple of the canonical class $K_S$ of our surface $S$.

To compute this invariant we resort to algebraic geometric techniques. The
appropriate moduli parametrize stable bundles together with sections. Actually
we also have to consider a second moduli parametrizing stable sheaves together
with maps to the canonical line bundle. These objects lead to
projective moduli spaces $MP^H(c_1,c_2)$.
The $\mu$--classes  on  $MP^H(c_1,c_2)$ are pulled back from
the corresponding $\mu$--classes on the moduli $\sM^H(c_1,c_2)$ of
stable sheaves along the natural map forgetting the section or cosection.
An intersection product of these algebraically defined classes gives the
same result as the differential geometric computation. This is a
consequence of the
corresponding statement for ordinary Donaldson polynomials (\cite{Mo},
\cite{Li}). The final result of our computation then is:

\proclaim {Theorem} For a simply connected Dolgachev surface $S$ with multiple
fibres of multiplicities $p$ and $q$ the coefficients $a(n)$ and $b(n)$ are
given by
$$\align
a(n)=&\,3n\\
b(n)=&\,(2p^2q^2-2p^2-2q^2-1)n.
\endalign
$$
\endproclaim

The coefficient $b(n)$ thus carries exactly the same information as
does Donaldson's $\Gamma$--invariant for Dolgachev surfaces (cf. \cite {Ba1}).
It suggests a strong relationship between the Spin--invariants and
Donaldson's original ones. Our methods are not good enough to get
the third coefficient, which we expect to finish the
diffeomorphism classification of Dolgachev surfaces. As we understand,
this classification was recently obtained by R. Friedman \cite{Fr} by different
methods.

Nevertheless the methods developped in this paper are of interest in their
own. The main technical problems in dealing with the moduli spaces are:
They have the wrong dimension and furthermore consist of various components
most of which are singular at every point. To compensate for the wrong
dimension
we use bundle theory. This replaces counting of intersection points with
determining an appropriate top Chern class. The nonreduced structure is
determined by explicitly constructing infinitesimal extensions and showing
that any further extension is obstructed. This explicitness is necessary in
order to understand whether and if so which sections or cosections extend.

The paper is organized as follows: The next chapter briefly reviews the
definition of the Spin--polynomials and provides machinery for dealing with
the singularities in the moduli spaces. The third chapter then gives the
definition of the diffeomorphism invariant for Dolgachev surfaces. Chapter four
does
the necessary bundle theory, which resembles the one developped in \cite{Ba2}.
The final
chapter then handles the combinatorics.
\bigskip\bigskip
\heading
II. $\Spin^{\Bbb C}$ polynomial invariants and wall structures
\endheading
\medskip
\subheading{IIa $\Spin^{\Bbb C}$ polynomial invariants}
\smallskip
Let $X$ be a four dimensional Riemannian manifold. An integral lift
$C$ of the second Stiefel-Whitney class $w_2 (X)$ defines a
$\Spin^{\Bbb C}$--structure on $X$, which in turn
defines a pair $W^+ , W^- $ of rank--2 complex Hermitian vector bundles
satisfying:
$ c_1 ( \Lambda ^2 W^{\pm} ) = C$.
A choice of a connection $\nabla$ on the line
bundle $\Lambda ^2 W^+ $ gives rise to a Dirac operator
$$ \Cal D^{C,\nabla} : \Gamma (W^+ ) \rightarrow \Gamma (W^- ), $$
which, when coupled with a connection $a$ on a Hermitian rank--2  bundle $E$,
gives a Fredholm operator
$$ \Cal D_a ^{C,\nabla} :
\Gamma (W^+ \otimes E ) \rightarrow \Gamma (W^- \otimes E).$$
We will assume its index to be nonpositive,
$$ \ind( \Cal D_a ^{C,\nabla} ) =
\dim \ker (\Cal D_a ^{C,\nabla}) - \dim \coker (\Cal D_a ^{C,\nabla}) \le 0 .$$
The family $\Cal D_a^{C,\nabla}$ of  Dirac operators over the space
$\Cal A $ of hermitian connections on $E$ can be used to define a subspace
of the moduli space of anti self dual (ASD) connections (see also [PT]):

\definition{2.1 Definition }
Fixing $C$, $\nabla$ and the metric $g$ on $X$, the
moduli of 1--instantons is the subspace
$$ \Cal M_1 ^{g,C,\nabla} =
\{ [a] \in \Cal M_{ASD}^g | \dim( \ker \Cal D_a ^{C,\nabla}) >0 \} $$
of the moduli of ASD connections on the adjoint $SO(3)$-bundle, which as
usual is identified with the moduli of Hermite-Einstein bundles on $E$.
\enddefinition

This moduli space provides invariants of the smooth structure of the manifold
$X$ (compare [PT] for an integer and [Ty] for a polynomial invariant).
It is convenient to also consider another moduli space,
which should be thought of as a resolution of those singularities of
$ \Cal M_1 ^{g,C,\nabla} $ which arise from a big kernel of of the
Dirac operator. This space has to be used in order to apply
general position  and bordism arguments to the moduli of 1--instantons
(cf. ch.1 of [PT]). The map which assigns to a connection $a$ and a
section $\sigma$ of the bundle $(W^+ \otimes E)$ the self dual part of the
curvature of $a$ and the evaluation of the Dirac operator on $\sigma$
$$ (a, \sigma) \mapsto (F_a ^+ , \Cal D_a^{C,\nabla} (\sigma) ) $$
is equivariant with respect to the standard action of the gauge group
$\Cal G_{E}$ of the bundle $E$. It therefore
gives rise to a section $s$ (see 1.1.26 of [PT])
of a corresponding (Hilbert) bundle over the space
$$ \Cal P = \Cal A_{\lambda} \times_{\Cal G_{E}}
S ( \Gamma (W^+ \otimes  E) ). $$

\definition{2.2 Definition }
The moduli space of pairs $ \Cal {MP}^{g,C,\nabla} $
is the zero set of this section.
A point in this moduli of pairs thus consists of a pair of objects: an
ASD connection and an element of the projectivisation of the kernel of
the coupled Dirac operator.
(Note that the central $U(1)$ in the gauge group $\Cal G_E$
acts trivially on the first factor in $\Cal P$ and in the standard
way on the second.)
\enddefinition

\proclaim{2.3 Transversality Theorem \cite{PT, ch 2.3}}
For a generic metric $g$ on $X$ and a generic connection $\nabla$ on
$\Lambda ^2 W^+ $ the moduli $ \Cal {MP}^{g,C,\nabla} $ of pairs is a smooth
manifold of dimension $$ \dim(\Cal {MP}^{g,C,\nabla}) =
2(4c_2(E) - c_1(E)^2) +2(\ind \sD_a^{C,\nabla} - 1).$$
\endproclaim

There is an obvious map
$ \pi :\Cal {MP}^{g,C,\nabla} \to  \Cal M_1 ^{g,C,\nabla} $
with fibre $\Bbb P (\ker \Cal D_a ^{C,\nabla} )$ over the
point $ [a] \in \Cal M_1 ^{g,C,\nabla}$. The tangent space at a point in the
moduli of pairs is conveniently described in the following exact sequence:
$$ 0 \to T_{a,\sigma} \Cal {MP}^{g,C,\nabla} @>>>
(\ker \Cal D_a^{C,\nabla})/<\sigma> \oplus H^1(AdE) @>i>>
\coker \Cal D_a^{C,\nabla} $$
with $ i(\sigma', \alpha) = \varpi(\alpha * \sigma) $.
Here $*$ is multiplication of (vector valued) spinors with (matrix valued)
1--forms, followed by spinor multiplication and
$\varpi$ is an orthogonal projection to the space of
negative harmonic spinors $\coker \Cal D_a$.
Finally we denote by $<\sigma>$ the line in a vector space
corresponding to a point $\sigma$ of its projectivisation. This exact sequence
globalizes to a fibre of $\pi$ as an exact
sequence of bundles over  $\Bbb P (\ker \Cal D_a ^{C,\nabla} )$ :
$$ 0 \to T \Cal {MP}^{g,C,\nabla}_{|\pi ^{-1}a} @>>>
 T \Bbb P (\ker \Cal D_a ^{C,\nabla} ) \oplus H^1(AdE) @>>>
\coker\sD_a^{C,\nabla}\otimes\sO_{\Bbb P(\ker\sD_a ^{C,\nabla})}(1).$$
A neighborhood of the point $[a]$ in the moduli space of ASD--instantons can
be described as a neighborhood of the origin in the zero set of
the obstruction map
$$ \psi : H^1_a (AdE) \to H^2_a (AdE) $$
in the Kuranishi description of the moduli space (cf. [FU], cor. 4.8).
Here $ H^i _a (AdE) $ denotes the homology in the Atiyah complex
 $$\Omega^0 (AdE) @>d_a>> \Omega^1 (AdE) @>d_a>> \Omega^2_+ (AdE) $$
associated to the connection operator $d_a$.
There is parallel description for the moduli space of pairs:

\proclaim{2.4 Lemma}
A neighborhood of the point $([a], \sigma)$ of the space
$\Cal {MP}^{g,C,\nabla}$  can be described as a neighborhood of
the zero set of the obstruction map:
$$  (\ker \Cal D_a)/<\sigma> \oplus H^1(AdE) @>\Psi>>
\coker \Cal D_a^{C,\nabla} \oplus H^2(AdE), $$
$$ \Psi(\rho, \alpha') = (\varpi(\alpha * \sigma) + O(|\alpha|^2) , \psi). $$
This globalizes to a neighborhood  of the
 $ \Bbb P ( \ker \Cal D _a)\times o\subset
 \Bbb P ( \ker \Cal D _a) \times H^1_a(AdE) $
which can be be described
through the zero set of an obstruction section of the bundle
$$ \coker \Cal D_a ^{C,\nabla} \otimes
\Cal O_{ \Bbb P (\ker \Cal D_a ^{C,\nabla} )} (1)
\oplus H^2_a(AdE). $$

\endproclaim

\demo{Proof}
 By [FrU], lemma 4.7, there is a diffeomorphism $\chi$ of
$ \Omega^1 (AdE)$ such that
$$ F_{a+\chi(\alpha)} ^+ = d_a ^+ (\alpha) + \psi(\alpha). $$
The linearisation of this diffeomorphism at zero is the identity map.
Furthermore,
$$ \Cal D_{a+\chi(\alpha)}^{C,\nabla} (\sigma + \rho) =
\Cal D_a^{C,\nabla} (\rho) + \chi(\alpha) * \sigma. $$
Extend the diffeomorphism $\chi$ to a diffeomorphism of the space
$ \Omega^1 (AdE) \times \Gamma (W^+ \otimes E ) $:
$$ (\alpha, \sigma') \mapsto (\chi(\alpha), \chi_{\alpha}(\sigma')) $$
to get the standard Kuranishi description:
$$ \Cal D_{a+\chi(\alpha)} (\sigma + \chi_{\alpha}(\sigma')) =
\varpi( \chi(\alpha) * \sigma) .$$
Now the leading term in the first component of $\Psi$ comes from the leading
term of the nonlinear term :
$$\varpi( \chi(\alpha) * \sigma) =
\varpi (\alpha * \sigma) + O(|\alpha|^2)* \sigma.\qed$$
\enddemo

The map $\Psi$ is not an honest obstruction map
since it has the linear term. However, it is convenient since it is explicit.

This lemma in paricular makes it possible to define an orientation of
the space $\Cal {MP}$ at a smooth point having fixed an orientation of
the moduli space of instantons. The latter is given by an orientation of
$H^1_a(AdE)$. Using the standard orientation of complex projective space one
gets an orientation of
$ \Bbb P ( \ker \Cal D _a) \times H^1_a(AdE) $. This induces a natural
orientation on our moduli space of pairs, as it is the zeroset of a
(transversal) section of a
complex vector bundle on an oriented manifold.
In the special case of a Hodge metric the orientation of the
moduli of ASD-connections will be chosen  in
such a way that it coincides with the one arising from
the complex structure on the moduli space of stable bundles.

In what follows we fix the Chern classes $c_i$ of the bundle and sometime
omit them from subscripts.
Denote by $ \Cal V_{\Sigma} ^{g,C,\nabla} $ the codimension-2 submanifold
of the manifold $ {\Cal M}_1  ^{g,C,\nabla}$ dual to $\mu( \Sigma)$.
It is a lift of a codim-2 submanifold in the moduli of bundles over
the Riemann surface $\Sigma \subset X$ (cf. [Do3]). Let
$ \Cal {VP}_{\Sigma} ^{g,C,\nabla} $ be its preimage in
$ \Cal {MP}^{g,C,\nabla}$. The submanifold  $\Cal V_{\Sigma} ^{g,C,\nabla} $
is obtained as the zero set of a certain section $ s_{\Sigma} $
of the line bundle $ \Cal L_{\Sigma_i}$ representing $\mu(\Sigma)$. Using
such submanifolds in the moduli of pairs
one can apply the yoga of the Donaldson invariants
to define the $\Spin^\Bbb C$  polynomial invariants (cf. [Ty]).
As usual, genericity assumptions on the metric have to be made.
In particular there should be no
reducible connection in the moduli space (i.e. $b_2 ^+ (X) \ge 1$).
Suppose the moduli of pairs $\Cal {MP}^{g,C,\nabla}$ is a manifold of
dimension $2d$.
The Spin polynomial invariant $q_{c_1,c_2}^C\in\Sym^d(H^2(x))$ then is
defined (cf. \cite{Ty}) by counting oriented intersection points
$$q_{c_1,c_2}^C(\Sigma_1,\dotso,\Sigma_d)=
\#\{\cap   \Cal {VP}_{\Sigma_i} ^{g,C,\nabla}\}.\tag{$*$}$$

\proclaim{2.5 Theorem \cite{Ty}} Suppose $b_2^+(X) > 1$. Then the
 Spin polynomial invariant
$q_{c_1,c_2}^C$ is independent of the metric $g$ or the connection
$\nabla$. For a diffeomorphism
$f: X \to Y$ of 4-manifolds one has
$$f^* q_{Y,c_1, c_2}^C = q_{X,f^*c_1, f^*c_2}^{f^*C}.$$
The dependence on the $\Spin^\bC$ structure is described in the
formula
$$
q_{c_1,c_2}^{C+2\delta}=q^C_{c_1+2\delta,c_2+c_1\delta+\delta^2}.
$$
\endproclaim

In case $b_2^+(X) = 1$ the right hand side of $(*)$
actually depends on the metric. For small values of
$p_1 = c_1 ^2 -4c_2$ this dependence can be controlled.
In particular this leads to chamber structures as in \cite{Do1}.
A special case will be considered below.
\medskip
\subheading{II.b The Case of a Hodge metric}
\smallskip
Consider now the case where our favorite 4-manifold  $X$ has the structure
of an algebraic surface $S\subset\Bbb P^N$. In this case the moduli of
irreducible
ASD connections can be identified with the moduli of stable holomorphic
bundles (cf. \cite {Do2, UY, DK}). The natural choice of a $\Spin^{\Bbb C}$
structure then is the anticanonical class $C = c_1 (S) = -K_S $.
With this choice the kernel of the coupled Dirac operator is isomorphic
to the sum of the zeroth and second cohomology group of the
corresponding stable bundle. The cokernel is isomorphic to its first
cohomology group.
So the existence of a nontrivial kernel for the Dirac operator
coupled with an ASD--connection is equivalent
to the nonvanishing of either the zeroth or
the second cohomology group of the corresponding stable bundle. This leads
to a description of the moduli of pairs as the moduli space $MP^H(c_1,c_2)$
 of stable
bundles together with a section or a cosection. By a cosection of an algebraic
bundle $\sE$  we mean a section of $\sE ^{\vee} \otimes K_S$. We will denote by
$MP_{a,b}^H(c_1,c_2)$
the stratum of $MP^H(c_1,c_2)$ where $h^0(S;\sE)=a$ and $h^2(S;\sE)=b$. It maps
to the stratum $M_{a+b}^H(c_1,c_2)$ of the moduli space $M^H(c_1,c_2)$
of stable bundles. For simplicity, we will always assume the first Chern
class $c_1$ to be odd and the very ample divisor
$H$ (or the Hodge metric $H$) to be suitably generic such that there are
no properly  slope semistable sheaves (or reducible connections) in
the compactified moduli. Also for simplicity assume that $MP_{1,1}$ is empty.
Then one can compare the orientations of these components
defined in the
previous section with their natural orientations as
complex spaces. It turns out that the
orientation of $MP_{1,0}$ will coincide with its natural complex
orientation and the orientation of $MP_{0,1}$ will differ from the
natural complex one by
$(-1)^{\chi (E) +1}$.

For a description of $MP_{1,0}^H$ in algebraic geometric terms we will make
simplifying assumptions: Suppose $c_1$ is odd, so there are no properly
slope semistable sheaves and the moduli space $M^H(c_1,c_2)$ of stable
sheaves is projective. Assume furthermore the existence of a universal sheaf
$\sE_M$ on $S\times M^H(c_1,c_2)$. Then $MP_{1,0}^H(c_1,c_2)$ and
$MP_{0,1}^H(c_1,c_2)$ are defined as projective cones
$$
\Proj (\sE xt^2_{p_M}(\sE_M,p^*_S\omega_S)) \,\,\,\,\,\,\text{and}
\,\,\,\,\,\,\,
\Proj (\sE xt^2_{p_M}(\sE_M,p^*_S\omega_S)).
$$
(We denote by $\bP(V)$ the space of lines and by $\Proj(V)$ the space of
hyperplanes, which agrees with the Grothendieck definition, cf. \cite{Ha,
II.7}.) Here we take $\Sym^.(\sE xt^2_{p_M}(\sE_M,p^*_S\omega_S))$ as
the graded sheaf with $\Sym^0(.)=\sO_{\supp(.)}$. The maps $p_S$ and $p_M$
denote the projection of the product $S\times M^H(c_1,c_2)$ to the respective
factors.

This definition leads
to a nice behaviour with respect to base change: For  $g:Y\to M^H(c_1,c_2)$
we have
$$g^*\sE xt^2_{p_M}(\sE_M,p^*_S\omega_S)\cong
\sE xt^2_{p_Y}(\sE_Y,p^*_S\omega_S)
$$
because the fiber $S$ of $p_M$ is two dimensional and smooth (cf.
\cite{BPS,2.2}).

A Hodge metric and a harmonic connection on $\Lambda^2 W^\pm$ in general are
not generic and thus the transversality theorem usually fails to hold. A method
to deal with such a situation was developed in \cite{PT, ch 3} in the special
case where no $\mu$--classes appeared. But the $\mu$--classes can be included
without any difficulties into this method by replacing the section $\tilde s$
of the bundle $\sF$ in \cite{PT, 3.1.17} by the section $\tilde s+\sum_i
s_{\Sigma_i}$ of $\sF\oplus(\oplus_i\sL_{\Sigma_i})$. Applying the same
arguments one gets the following analog of \cite{PT, prop. 3.2.4}:

\proclaim{2.6 Proposition}
 Denote by $\mu^H$ the intersection of  $\Cal{VP}_{\Sigma_i}$, $i=1,\dotso, d$,
for the Hodge metric in consideration, i.e.
$\mu^H \subset MP_{1,0} \cup MP_{0,1} $. Set $\mu^H_{1,0}=\mu^H\cap MP_{1,0}$.
For a generic deformation of the parameters $g$ and $\nabla$ the following
formula gives the ( 0-dimensional) fundamental class of the manifold
$\mu ^{g,C,\nabla}_{1,0}$ which is the deformation of $\mu^H_{1,0}$:
$$ [\mu ^{g,C,\nabla}_{1, 0}] = [c_* (\ind At \oplus \ind \Cal D \oplus (\oplus
\Cal L _{\Sigma_i}) _{|\mu^H _{1,0}}) \cap (\can(\mu^H_{1,0}))]_0.\qed$$
\endproclaim

Here, $c_*$ denotes the total Chern class and $\can(Y)$ the canonical class
(\cite{Fu, Ex. 4.2.6}) of a variety $Y$.
An analog statement is true for a deformation of
$\mu ^{g,C,\nabla}_{0,1} = \mu ^{g,C,\nabla} \cap MP_{0,1}$.
Instead of using the stable bundle $E$ one has to use $E^* \otimes K_S$.
When combining the two contributions, one as of course to take the
orientations into account.

In the example of Dolgachev surfaces below the stratum
$\mu^H_{1,1}$ in general will not be empty, but $\mu^H_3$ will be.
(This happens for small $\vert n\vert<pq$ only).
That $\mu^H_{1,1}\setminus\mu_3^H$ actually deforms into the empty set
if the orientations do not conincide, follows from the following consideration:
The local model of a point of degree $n$ in $\mu^H_1$ in this case will be
$V \to Hom( H^0(E) \oplus H^2(E), H^1(E))$
given by the matrix
$ \left( \matrix z^n & 0 \\ 0 & \overline z ^n \endmatrix \right) $
(cf. [PT], ch.1, formulae 1.5.5 - 1.5.8) where $z$ is a complex
variable in the 1-dimensional tangent space $V$ to $\mu^H$.
This local model
can be easily deformed to a family of matrices without degeneration. Simply
take
$ \left( \matrix z^n & i \epsilon \\ i \epsilon & \overline z ^n \endmatrix
\right) $
with a nonzero $ \epsilon $.

There is one particular case of interest, where the moduli space $\sM_1$ of
1--instantons has the expected dimension, but the rank of the kernel of the
coupled Dirac operator is greater than one. Assume again for simplicity
that $\sM_{1,1}$ is empty. We want to compare the intersection of
$\mu$--classes in $ PM_{1,0}^H$ with the corresponding  intersections of
$\mu$--classes in $\sM_{1,0}^H$. The result roughly states that a point in
$\sM_{1,0}$, which is represented by a bundle $\sE$ gives $h^0(\sE)$ points
in $MP_{1,0}$.

\proclaim{2.7 Corollary } Suppose $\sM_1$ is of expected dimension $d$
and let $Z\subset \sM_{1,0}$ be an irreducible component.
Let $\sE_\zeta$ be the pull back of a universal sheaf on
(a Zariski open subset of)
$Z$ to its generic point $\zeta$ and let $h$ denote the length
of the relative $\Ext$ sheaf $\sE xt^2_{p_\zeta}(\sE_\zeta,p^*_S\omega_S)$.
If $k$ is the
algebraic intersection number of Cartier divisors $\sV^H_{\Sigma_i}$,
$i=1,\dotso,d$, in $Z_{red}$, then the algebraic intersection number of
the Cartier divisors $\Cal{VP}^H_{\Sigma_i}$ is $hk$.
\endproclaim
\demo{Proof} Since the moduli spaces are projective, we may assume the divisors
$\sV^H_{\Sigma_i}$
to be very ample and in general position, so that the finitely many
geometric points of the intersection subschemes
are all contained in the Zariski open subset of $Z_{red}$
 parametrizing sheaves in with the lowest possible number of sections.
We can now apply the local description of $\mu^H$ as the
zero locus of a section of a bundle
$$ F =   (\oplus \Cal L _{\Sigma_i}) \oplus H^1(\sE)\otimes\Cal O_{ \Bbb P
(H^0(\sE))} (1)  \oplus H^2(\sE nd_0(\sE)) $$
on
$Y= \Bbb P ( H^0(\sE)) \times H^1(\sE nd(\sE))$. The equality
$$ (F - TY)_{|\mu^H} = (\ind At \oplus \ind \Cal D \oplus (\oplus \Cal L
_{\Sigma_i}) _{|\mu^H})$$
implies that the total Chern class of this virtual bundle is $1$.
So we need only determine $[can(\mu^H)]_0$. If $Z$ and thus the point in
$\sM_{1,0}^H$ is reduced, this canonical class is just the Euler characteristic
of the projective space $\bP (H^0(\sE))$. Otherwise we can analyse the
scheme $\mu^H$ locally over a nonreduced point $\xi$ in $\sM_{1,0}$, where
in general there are imbedded subschemes: Inductively, over subschemes $\xi_n$
of $\xi$, we see that the obstruction to extending $\pi^{-1}(\xi_n)$ to a
scheme
over $\xi_{n+1}$ vanishes along a projective subspace of $\bP
(H^0(\sE_{\xi_n})^\vee)$.
The canonical classes add.
\qed
\enddemo

\medskip
\subheading{II.c Wall structure}
\smallskip
In the case $b_2 ^+(X) = 1$ reducible ASD-connections will appear in
one-parameter families of metrics on $X$ although they are absent
in the moduli space for generic metrics. So there will be a
singularity in the bordism moduli space defined by a path of metrics
if the path meets a metric
which allows for reducible ASD-connections (cf. [D1]).
Any decomposition $ E = M \oplus (c_1 - M)$
with  $c_1^2-4 M \cdot (c_1 - M) \ge p_1 (adE)$
defines
a one-codimensional subspace in the space $\sC$ of all metrices.
Sending a metric to a selfdual harmonic form gives a period map
$\sC\to\bP(H_+^2)$ to the projectivisation of the cone $H_+^2\subset
H^2(X;\bR)$ of positive vectors. A decomposition of $E$ as above defines a wall
$(c_1-2M)^\perp$ in $\bP(H_+^2)$.  A connected component in the complement
of the set of all walls is called a chamber.
In the case of moduli space of pairs not every
wall is effective in the sense that intersections of $\mu$--classes give
different results on both sides of the wall.

Our aim here is to prove that for $p_1 = c_1 ^2 -4c_2 > -7$ or
$p_1 = -8, w_2 \ne 0$ the Spin--polynomial invariants depend only on the
chamber to which the metric maps.

Consider a reducible ASD-connection of type $ M \oplus (c_1 - M) $ with
$(c_1 -2 M)^2 = p_1 (adE)$, $c_1\cong w_2(adE)\mod 2$. It lies
in the moduli space itself, not in its Taubes compactification.
If a path of metrics crosses the  wall $(c_1-2M)^\perp$
then the corresponding instanton bordism
$ \cup_{t} \sM_{ASD}^{g_t} $ will have a singularity at the point given by
$ M \oplus (c_{1} - M) $. A neighborhood of this singularity in the moduli of
ASD-connections can be described as a (real) cone over the
complex projective space
$ \Bbb P (H^1 (Ad( M \oplus (c_{1} - M) )) $  (without loss of generality
we can assume the second cohomology to vanish). In other terms we have a smooth
bordism between the moduli of ASD-connections
$$ {\Cal M}_{ASD}^{g_{-1}} \,\,\,\,\,\,\, and {~~~}\,\,\,\,\,\,\,\,
 {\Cal M}_{ASD}^{g_1} \cup
\bP(H^1 (Ad( M \oplus (c_1 - M ))) .$$

The same picture is true for the moduli spaces of 1--instantons with only a
minor modification: Now the 1--instantons describe a homology class
of $\bP(H^1 (Ad( M \oplus (c_1 - M) ))$. If this homology class happens to
vanish, then the wall has no effect on the $\Spin ^{\Bbb C}$--invariants.
In order that this homology class does not vanish one has to have  either
$ \ind {\Cal D}_{2M} > 0 $ or $ \ind {\Cal D}_{2(c_1 - M)} > 0 $ (cf.
\cite{PT,1.4.2}).
Recall that the index of the Dirac operator on a manifold $X$ with
$\Spin ^{\Bbb C} $--structure $C$ coupled with a connection on a line
bundle $L$ is given by the formula
$$ \ind \Cal D _L = \left((C+c_1(L))^2-\Sign(X)\right)/8 .$$

\noindent{\bf 2.8 Remark.}
This condition is effective and will be applied lateron
to a Dolgachev surface $S$. We will take
$c_2 = 2$ and $c_1$ a rational multiple of the canonical class $K_S$.
Suppose $M$ is perpendicular to $K_S$. Then one gets
$$ \ind \Cal D _{2M} = \left((C+2M)^2-\Sign(S)\right)/8 = 4M^2-8=0$$
and similarly $ \ind \Cal D _{2(c_1 -M)} = 0$.
So these particular walls have no effect on the $\Spin^\bC$--invariants.

To analyse the link of the singularity in the bordism moduli space
one can apply Porteus formula \cite{Fu, ch. 14}.
The link comes as the  fundamental class
of the degeneracy locus  $$ D = {\Cal M}_1^C \cap
P(H^1 (Ad( M \oplus (c_{1} - M) )) $$ of a morphism between bundles
$ E \to F $:
Consider the index bundle of
the family of coupled Dirac
operators over the link of the singularity in the bordism moduli space. This
index bundle can be described as an equivariant complex of $ S^1 $--bundles
$$ \ker {\Cal D} _{ ( M \oplus (c_1 - M) )} \rightarrow \coker
{\Cal D} _{ ( M \oplus (c_1 - M) )} $$
(with linear actions of weigths $+1$ and $-1$ on the former and the latter)
over the $ S^1 $--space $ H^1 (Ad( M \oplus (c_1 - M) ) $ (with linear action
of weight two). Because of the
transversality theorem \cite{PT, 1.3} we may assume:
$$ \ker {\Cal D} _{ (c_1 - M) }= \coker {\Cal D} _{  M } = 0 .$$
We twist the above complex by a one-dimensional $S^1$--representation of weight
$1$ to get representations which descend to bundles over
$\bP H^1(Ad(M\oplus (c_1-M)))$.
These bundles can be written as a formal difference
$$ \left(\ind({\Cal D} _{(c_1 - M)} ) \otimes {\Cal O }
(-1)\right)-\ind({\Cal D} _{M}) .$$
and Porteus formula gives for $k = \codim {~} {\Cal M}_1^C $
$$ [D] = c_k ( \ind({\Cal D} _{(c_1 - M)} ) \otimes {\Cal O } (-1) {~} -
{~} ( \ind({\Cal D} _{M}))) {~}.$$

This fundamental class is given in terms of
indexes and in particular is independent of
any metric. To get formulae describing the change in the polynomial invariant,
one of course has to evaluate products of $\mu$--classes on this homology
class.
However, the result again is independent of any metric. It only depends
on the wall defined by the reducible ASD-connection.

This method fortunately works in one more case to the effect
that the invariants for $p_1(adE) = -8$ are functions depending
only on chambers in cohomology: Suppose now the reducible connection $A_r$ on
the bundle $ M \oplus (c_1 - M) $ to satisfy
$ M.(c_1 - M) = p_1 (E) +4 $. This means
it lies in the top stratum of the compactification of the moduli space.
First we describe the link of this singularity in the compactification.

The $ \Spin^\bC (4) $--structure on our
manifold determines a principal $U(2)$-bundle $P^+$ corresponding to
$W^+$. Denote by $P_{M \oplus (c_1 -M)}$ the principal bundle corresponding
to our reducible bundle. Reducibility provides an isomorphism
$$ P_{M \oplus (c_1 -M)} = U(2) \times _{T^2} (P_M \times P_{(c_1-M)}),$$
where $ P_M, P_{(c_1-M)} $ are the principal bundles corresponding to
$M$ and $c_1-M$, respectively.
  Let $ Gl $ be the subspace of all hermitian isomorphisms in the bundle of
homomorphisms
 $ Hom (M \oplus (c_1 - M) , W^+) $ (recall
both arguments are equipped with hermitian metrices). The space Gl is a
principal $T^2$--bundle
$$ (P^+ \times P_{M\oplus (c_1-M)})/U(2) =
(P^+ \times (P_M \times P_{(c_1-M)}))/(T^2) $$
over the space $P^+ /(S^1 \times S^1) = \Bbb P (W^+).$
Here $ T^2 \subset U(2) $  is the holonomy centralizer of the reducible
connection. Let $ S^1 \subset U(2) $ be the center, considered as
 holonomy centralizer
of a 1-instanton over $ S^4 $. It also acts on the space $Gl$.
The "diagonal" subgroup
$ S^1 _{\triangle} \subset T^2 \times S^1 $, defined by
$$ e^{i \phi } \mapsto \left( {~~} \left( \matrix  e^{i \phi }
& 0 \\ 0 &  e^{i \phi } \endmatrix  \right)    {~~} ,
 {~~}   \left( \matrix e^{-i \phi } & 0 \\ 0 &  e^{-i \phi }
\endmatrix   {~~}  \right)  \right) $$
acts trivially on $Gl$ and in fact we have an action of
$$ T^2 = T^2 \times S^1 / S^1 .$$

We will use the space $Gl \times \bR_{\ge0}$ (where $ \bR_{\ge0}$
parametrises the concentration of a glued instanton
on the 4-sphere) as the space of parameters needed to describe the
neighborhood of ideal connections in the Taubes compactification
(cf. \cite{Do3, ch. 5}). Since in the Taubes compactification all
glueing parameters disappear for ideal connections,
 one has to consider the space
$ Gl \times \bR_{\ge0} / Gl \times {0} .$
The neighborhood of an ideal connection of type $(A_r, x)$
for $x \in X$
in the compactified moduli space is described by
$$ (Gl \times \Bbb R_{\ge0} / Gl \times {0}) \times H^1(ad(A_r)) / T^2 $$
(cf. Th 4.53 in \cite{loc. cit.}).

\proclaim{2.9 Lemma} The space
$$ (Gl \times \Bbb R_{\ge0} / Gl \times {0}) \times H^1(ad(A_r)) / T^2 $$
is a fibration $\Bbb V$ over $X$ with fiber a real cone over an algebraic
subvariety $V$ of $\Bbb P (S^2(E^*)) \otimes H^1(ad(A_r))) $.
The variety $V$ is the image under the Segre embedding
$$\Bbb P (S^2(E^*)) \times \Bbb P(H^1(ad(A_r)) \to
\Bbb P (S^2(E^*) \otimes H^1(ad(A_r))$$
of the subspace
$R \times \Bbb P(H^1(ad(A_r))$
where $R$ is a rational
norm curve in $\Bbb P ( S^2(E^*)) $.
\endproclaim
\demo{Proof} First note that the diagonal $S^1 \in T^2$ acts trivially on
$H^1(ad(A_r))$, so we can factor it out to get
$$ ((Gl/S^1) \times \bR_{\ge0} / (Gl/S^1) \times {0}) \times H^1(ad(A_r)) / S^1
.$$
The space $Gl/S^1$ is isomorphic to the principal $S^1$-bundle associated to
the line bundle $T_p \Bbb P(W^+) \otimes p^*\sL^{-1}$. Indeed, the principal
$S^1$-bundle $P^+/S^1$ is associated with the line bundle $T_p \Bbb P(W^+)$ on
$\Bbb P(W^+)$ and the principal $S^1$-bundle $ (P_M \times P_{(c_1-M)}))/(S^1)
$ is
associated with a line bundle $\sL$ with Chern class $2M-c_1$
on $X$. Therefore
$$((P^+ \times (P_M \times P_{(c_1-M)}))/(T^2))/S^1 =
((P^+/S^1) \times (P_M \times P_{(c_1-M)})/S^1)/(S^1) $$
and the statement follows from the definition of tensor product of line
bundles.
Now one can identify the space
$$ ((Gl/S^1) \times \bR_{\ge0} / (Gl/S^1) \times {0}) $$
in the following way. Take the
obvious $S^1$-equivariant map
$$f: ((Gl/S^1) \times \bR_{\ge0}  \to T_p \Bbb P(W^+) \otimes p^*\sL^{-1}.$$
It projects $ (Gl/S^1) \times {0}) $ to $ Gl /T^2 = \Bbb P (W^+)$.
We now construct a $S^1$-equivariant map of total spaces
$$v: T_p \Bbb P(W^+) \otimes p^*\sL^{-1} \to  S^2(E^*) \otimes \det E \otimes
p^*\sL^{-1}$$
 as the twist of the obvious map of total spaces
$ \Cal O _{\Bbb P (E)}(2) \to S^2(E^*)$
with the line bundle $\det E \otimes p^*\sL^{-1}$.
It is easy to see that the composite $S^1$-equivariant
map $vf$ projects $ (Gl/S^1) \times {0}) $ to the zero-section of
$S^2(E^*) \otimes \det  E \otimes p^*M^{-2}$.

To include the vector space $H^1(ad(A_r))$ into the picture consider
the composite $vf \times id_{H^1(ad(A_r))}$:
$$ ((Gl/S^1) \times \bR_{\ge0} \times H^1(ad(A_r)) \to S^2(E^*) \otimes \det  E
\otimes p^*\sL^{-1} \times H^1(ad(A_r))$$
which leads to the Segre imbedding
$$ S^2(E^*) \otimes \det  E \otimes p^*\sL^{-1} \times H^1(ad(A_r)) \to
 S^2(E^*) \otimes \det  E \otimes p^*\sL^{-1} \otimes H^1(ad(A_r)) .$$\qed
\enddemo

\proclaim{2.10 Corollary} The singularity in the bordism of instanton spaces
$ \cup_{t} {\Cal M}_{1}^{g_t} $ given by a reducible connection of the
considered type is the space $\Bbb V$.
\endproclaim

\proclaim{2.11 Lemma } Suppose the characteristic classes of
the underlying bundle $E$
satisfy  $p_1(Ad\, E) = c_1(E)^2 -4c_2(E) > -7$ or
$p_1(Ad\,E) = -8, w_2(Ad\,E) \ne 0$. Then
 the $\Spin^\bC$--invariants depend only
on chambers, not on metrics.
\endproclaim
\demo{Proof} We have determined the links of the singularities in the bordism
moduli spaces. Now consider a path of metrics such that  the periods of
the metrics cross a wall. The change in the invariants then is computed
by evaluating certain cohomology classes on these links which
do not depend on the metric in which we cross the wall. This gives a unique
value for the Spin--polynomial for any metric with a given period.\qed
\enddemo
\bigskip\bigskip
\heading
A diffeomorphism invariant for Dolgachev surfaces
\endheading
We want to evaluate some $\Spin^\Bbb C$--invariant in the case of simply
connected Dolgachev surfaces, i.e. of relatively minimal elliptic surfaces
with geometric genus $p_g=0$. As was shown in the original paper of Dolgachev
\cite{Dol}, there are exactly two multiple fibers of coprime multiplicity
$p$ and $q$. Here we allow $p$ or $q$ to be 1, noting that in this case
the surface is rational: A canonical divisor is linearly equivalent to
$F-F_p-F_q$, where $F$ and $F_p$, $F_q$ denote a generic and the respective
reduced multiple fibers. The ray $\Bbb Q[F]\subset H^2(S;\Bbb Q)$ of rational
multiples of the Poincar\'e dual $[F]$ of $F$ is
generated by a primitive element
$k=1/pq[F]=a[F_p]+b[F_q]$ with integers satisfying $aq+bp=1$.

Let $q_S(n)$ denote the $\Spin^\Bbb C$--invariant
$q_{K_S+2nk,2}^{H}$, evaluated on a chamber in
$H\in H^2_+(S;\Bbb R)$ containing $k$ in its closure.
Recall that the invariant is a function on the
set of chambers, i.e. the components of
$$H^2_+(S;\Bbb R)\setminus\cup_Ie^\perp,\phantom{@>>>@>>>@>>>}
I=\{e\in H^2(S;\Bbb Z)\,\vert\, -8\le (c_1-2e)^2\le -1\}.$$
In the present case we can, however, say more:

\proclaim{3.1 Theorem} The polynomial $q_S(n)\in\Sym^4(H^2(S))$ is a
diffeomorphism
invariant of Dolgachev surfaces, i.e. if $F:S\to S'$ is a diffeomorphism of
Dolgachev surfaces,
then
$$f^*q_{S'}(n)(x_1,\dotso, x_4)=\pm q_S(n)(f_*x_1,\dotso, f_*x_4).$$
Furthermore, $q_S(n)$ is a polynomial of the form
$$q_S(n) = a(n)Q^2+b(n)Qk^2 + c(n)k^4$$
in the intersection form $Q\in \Sym^2(H^2(S))$ and in $k$.
\endproclaim

\demo{Proof} The ray $\Bbb Q[F]$ is invariant under diffeomorphisms for
nonrational Dolgachev surfaces, i.e. $f^*\Bbb Q[F]=\Bbb Q[F']$ for a
diffeomorphism between rational Dolgachev surfaces (\cite{FM1, thm. 6A}).
So the set of chambers containing $k$ in its closure is mapped to itself. It
remains to show that the invariant takes the same value on neighboring such
chambers. To see this note that $k^\perp\subset H^2(S;\Bbb Z)$ carries an even
form. So the respective chambers are separated by walls $e^\perp$ with
$e^2=-2$. The result now follows from the fact (cf. 2.8)
that chambers separated by
walls perpendicular to (-2)--classes $e\in k^\perp$ have the same invariant.
Thus $q_S(n)$ is an element in $\Sym^4(H^2(S))$ which is invariant under
the diffeomorphism group of $S$. It is clear that $\pm k$ is an invariant
linear form on $H_2(S)$. So it suffices to determine the invariant polynomials
on $\ker k$. Now $\ker k$ is isomorphic to the lattice $-\w E_8$ with
radical generated by $1/pq F$. By \cite{FM1, III.2} the transvection
$T_y(x)=x+(x.y)F$ can be realized by a diffeomorphism of $S$ for
$x,y\in -\w E_8$. From this it follows that invariant polynomial functions
on $-\w E_8$ are pulled back from invariant poynomial functions on
$-\w E_8/\rad=-E_8$. These are well known to be generated by invariant
polynomials in degree 2 (corresponding to $Q$), 8, 12, 14, 18, 20, 24, 30.
The theorem thus is established for nonrational surface. In the case of a
rational minimally elliptic surface one has to use the fact that the
diffeomorphism group acts transitively on the set
$$\{\kappa\in H^2(S;\Bbb Z)\, |\, \kappa^2=0, \kappa\, \,\text{indivisible},
\kappa^\perp\, \text{carries an even form}\}.$$
(cf. \cite{FM1, II.2.4} and \cite{Wa}.)
\qed\enddemo
\bigskip\bigskip
\heading
 Bundle theory, multiplicities
\endheading
\medskip
\subheading{IV.a Description of the bundles}
\smallskip
In this section the bundles relevant for the computation of the spin polynomial
invariant will be described. The discussion is completely analogous to the one
in \cite{Ba2, ch. II}, so we will be brief in our presentation.\par

Let $S$ denote a relatively minimal simply connected Dolgachev surface.
The surface  $S$ is
projective and we will assume it to be suitably generic. It will
suffice that $S$ be nodal in the sense of \cite{FM2, ch.2}: The multiple fibres
have smooth reduction and all singular nonmultiple fibres are irreducible
rational curves with a single ordinary double point.
The attribute {\it vertical} will be used for any divisor
linearly equivalent to a rational
multiple of a generic fiber. Such a vertical divisor is linearly equivalent
to  $lF+mF_p+nF_q$ for uniquely determined  integers
$l$, $m$ and $n$ with $0\le m<p$ and $0\le n<q$.

\proclaim{4.1 Lemma} The vertical line bundle $
\sO_S(lF+mF_p+nF_q)=\sL$ has cohomology groups of the
following dimensions:
$$h^0(S,\sL)=\max\{l+1,0\}, \qquad h^2(S,\sL)=
\max\{-l, 0\},$$ $$h^1(S,\sL)=\max\{l, -1-l, 0\}.$$
\endproclaim

\demo{Proof} By definition, $h^0(S,\sL)=h^0(\bP^1,\pi_*\sL)
=h^0(\bP^1, \sO_{\bP^1}(l))$. For $h^2$ use Serre duality, for $h^1$ Riemann
Roch.
\qed
\enddemo

Let $H$ be an ample divisor "near $K_S$" (i.e. assume
$H-2(H.F)F $ is ample, too).
We will be interested in those $H$--stable bundles $\sE$
with Chern classes $c_1=D$ a vertical divisor and $c_2=2$ which admit a section
or a cosection.
 These bundles
represent the closed points in the subscheme $\sM_1^H(D,2)$ of the
moduli space $\sM^H(D,2)$.
The choice of ample divisor implies that $\sE\in\sM_1^H(D,2)$
arises as an extension
$$0\ra\sO_S(C)\lra\sE\lra\sJ_Z(D-C)\ra 0\tag{${*}$}$$
with $Z\subset S$ a 2--cluster, i.e. a zero dimensional subscheme of
length $2$, and $C$ an effective
vertical divisor. In particular the direct
image sheaf $\pi_*(\sE)$ is locally free on $\bP^1$ of rank either 1 or 2.
The bundle will be called of typ 1 or 2, correspondingly.
These bundles have been described in \cite{Ba2, prop 3.4 and 3.6}:

\proclaim{4.2  Proposition (Properties of type 1 bundles)}
Let $\sE\in \sM_1^H(D,2)$ be a type 1 stable
bundle.
Then the divisor $C$ and the subscheme $Z$ are uniquely determined.
On the other hand, for a given
vertical divisor $C$ there exists a semistable bundle as
above if and only if $2C-D$ is linearly equivalent to $-F+\alpha F_p+\beta F_q$
with nonnegative integers $\alpha $ and $\beta$ such that $\alpha/p+\beta/q<1$.
In particular these bundles are paramerized by a Zariski--open subset
of $\Hilb^2(S)$, which has the expected dimension $4$.
The Zariski tangent space of the moduli space $\sM^H(D,2)$ at $\sE$ has
dimension $5+h^2(S;\sO_S(D-2C))$.\qed
\endproclaim

\proclaim{4.3 Proposition (Properties of type 2 bundles)}
Let $\sE\in \sM^H(D,2)$ be a type 2 stable
bundle. Then it arises as an extension
$$
0\ra \sO_S(C_1)\oplus\sO_S(C_2)\lra\sE\lra
   \sJ_{Z\subset I}(D-C_1)\ra 0.
\tag{${*}{*}$}
$$
Here $I$ denotes the minimal vertical divisor containing $Z$ and
$\sJ_{Z\subset I}$ denotes the corresponding ideal sheaf. The divisor
$I$ and the line bundles $\sO_S(C_i))$ are uniquely determined by $\sE$.
Furthermore, if the divisor $I$ contains a nonmultiple fiber $F$,
then $h^0(\sO_{F\cap Z})\ge 2$. The bundles $({*}{*})$ admit at most
4--dimensional moduli.\qed
\endproclaim

\demo{Remarks}
1. The subspace $\sM_1^H(D,2)$ of the moduli space $\sM^H(D,2)$ consists of
several components in general. The type--1 components are indexed by the pairs
$(\alpha, \beta)$ in the proposition. A listing of the different type--2
components can be obtained by considering which possible $Z\subset I$ may
occur. For a generic
bundles in these components $I$ is a generic fibre and the possible
$\sO_S(C_i)$ are characterized by the conditions
$H.C_1<H.(C_2+F)<H.(C_2+2F)$ and $h^0(\sO_S(C_1-C_2))=h^0(\sO_S(C_2-C_1))=0$.
These conditions are equivalent with: $C_1=C_2-F+\alpha F_p+\beta F_q$ with
$0<\alpha,\beta$ and $\alpha q+\beta p<pq$.
It turns out that the other strata of type--2 bundles are indeed contained as
subspaces in at least one of the components mentioned above.

It should also be pointed out that the singular locus of $\sM^H(D,2)$ contains
only type--1 and type--2 bundles. Moreover, the components of $\sM_1^H(D,2)$
as Weil divisors in $\sM^H(D,2)$ are not reduced in general. These singular
structures will be explored lateron.
\enddemo
\medskip
\subheading{IV.b Some families of stable sheaves}
\smallskip
We will evaluate the invariant $q_S(n)$ at four classes $x_1,\dotso, x_4\in
H_2(S;\bZ)$ with $x_4=F$. So it suffices to construct the divisor
$\mu(F)\subset
\sM_1^H(D,2)$ and the restriction of a universal sheaf to $S\times \mu(F)$.
For this we use the following observation in \cite{FM2, ch 3.9}: the
divisor $\mu(F)$ geometrically is represented by the locus
$$\mu (F)=\{[\sE]\in\sM_1^H(D,c)\,|\,h^0(\sE|_F\otimes\theta_F)\not=0\}$$
with $\theta_F\in\Pic^0(F)\setminus\{0\}$. In the present case, where the
bundles are all given as extensions
$$0\ra\sO_S(C)\lra\sE\lra\sJ_Z(D-C)\ra 0,$$
this means $\mu(F)=\{[\sE]\,|\,Z\cap F\not=\emptyset\}$. A component $C$ of
$\mu(F)$ then has to be counted with multiplicity
$\min_Ch^1(\sE|_F\otimes\theta_F)$, so
if $C$ generically parametrizes type--$n$ bundles, it has to be counted with
multiplicity $n$.

Here is the construction of a family parametrizing stable bundles which
generically are of the form
$$0\ra\sO_S\oplus\sO_S(-F+\alpha F_p+\beta F_q)\lra \sE\lra\sJ_{Z\subset F}\ra
0
$$
for $0<\alpha, \beta$ with $pq>\alpha q+\beta p$. Let $T$ denote the
tautological divisor in $F\times\Hilb^2(F)$. Then (compare \cite{Ba2, 3.7})
there is a unique nontrivial extension, which indeed is locally free:
$$
0 @>>> \sO_{S\times\Hilb^2(F)} @>>> \w\sE @>>>\sJ_T(\beta F_q, -K_{\Hilb^2(F)})
@>>> 0.
$$
Let $\tau:\sM_F\lra\Hilb^2(F)$ be the $\Bbb P^1$--bundle
$\Proj(pr_{2{*}}(\w\sE\vert_{\alpha F_p\times\Hilb^2(F)}))$. The family of
type--2 bundles then is described as follows (cf. \cite{Ba2,3.8}):

\proclaim{4.4 Proposition} Assume by symmetry $\alpha q<\beta p$. Then the
 kernel $\sU$ of the natural surjection
$$
(id\times \tau)^*\w\sE(\alpha F_p,0)   @>>>
 \sO_{\alpha F_p}(\alpha F_p)\boxtimes\sO_{\sM_F}(1)
$$
is a bundle parametrizing stable type 2 bundles.
The generic member of this family has a presentation
$$
0 @>>> \sO_S\oplus \sO_S(-F+\alpha F_p+\beta F_q) @>>> \sE
  @>>>\sJ_{Z\subset F} @>>> 0
.$$
The induced map of $\sM_F$ to the moduli space $\sM^H(\alpha F_p+\beta F_q,2)$
is injective on geometric points. The Chern classes of $\sU$ are
$$
\align
c_1(\sU) &=
 (\alpha F_p+\beta F_q)\times \sM_F -S\times \tau^*K_{\Hilb^2(F)}\\
c_2(\sU) &=\tau^*T + \alpha F_p\times c_1(\sO_\sM(1)).\qed
\endalign
$$
\endproclaim

Note that this construction was not symmetric in $p$ and $q$. This is reflected
in the computation (cf. \cite{Ba2, 3.10}) of products of the classes
$$
\mu(A)=A\backslash (c_2(\sU)-{1\over 4}c_1^2(\sU)).
$$

\proclaim{4.5 Corollary } Let $\sU$ be the bundle in (4.4).
Then for an effective divisor $A$ in $S$ we get
$$
\mu(A)^3=3\phi_2(\alpha,\beta)(A.F)(A.k)^2,
$$
where $k\in H^2(S;\Bbb Z)$ denotes the primitive class ${1\over pq}F$
and $\phi_2(\alpha,\beta)=\alpha(q-\beta)pq$.\qed
\endproclaim

The family $\sE_X$
of type--1 stable sheaves to be constructed will be over the blow
up $X$ of $S\times F$ along the diagonally imbedded curve $F\subset F\times F
\subset S\times F$ with exceptional divisor $E$. The sheaf will at a generic
point $x$
parametrize stable bundles of type 1 of the form
$$
0 \ra\sO_S(D)\lra\sE_x\lra\sJ_Z\ra 0
$$
with $D=-F+\alpha F_p+\beta F_q$ and $Z\in \Hilb^2(S)$, the image of $x$ under
the canonical map $X\ra \Hilb^2(S)$.

Set $d=pq-\alpha q-\beta p$ and denote by $\delta_p$ the integer
$$\delta_p=\left[\frac{d}{2q}\right]-1.$$
In \cite{Ba2, 3.16} a subscheme $\sR(\sL_p)$ of $pF_p\times pF_p$ was
constructed containing the diagonal subscheme
$\Delta_{2\delta_p}\subset2\delta_p F_p\times
2\delta_p F_p$. Let $\sJ_{\Delta\vert\sR(\sL_p)}$ denote the corresponding
ideal sheaf. The main properties of this subscheme are:
\proclaim{4.6 Lemma}
The space $\sR(\sL_p)$ is flat and finite of degree $\delta_p(\delta_p+1)$ over
$F_p\times F_p$ and
$$
\multline
\sE xt^*_{S\times X}(\sJ_{\Delta\vert\sR(\sL_p)}\boxtimes\sO_F,\sO_{S\times
X})=\\
\cases
\sO_{R(\sL_p)\times F}(-D,pr_1^*(D+(2\delta_p+1)F_p)& \text{for $*=2$}\\
0 &\text{else.}\qed
\endcases
\endmultline
$$
\endproclaim
This subscheme can be used to construct a geometric representativ of the second
Chern
class of $\sE_X$ (cf. \cite{Ba2, 3.18}):

\proclaim{4.7 Proposition} There is a unique nontrivial extension
$$0@>>>\sO_{S\times X}(D,L)@>>>\sE_X@>>>\sJ_{\Gamma}@>>>0$$
 giving
a flat $X$--family of stable sheaves. The subscheme $\Gamma\subset S\times X$
is the union of $(\sR(\sL_p)\cup\sR(\sL_q))\times F$ with
the graphs $\gamma(X)$ and $\rho(X)$ of the maps $\gamma,\rho:X\ra S$ factoring
through the projection maps $pr_i:S\times F\ra S$ ($i=1,2$). The divisor $L$
is given by $-L=E+\gamma^*(D+2\delta_p F_p+2\delta_q F_q-K_S)$.\qed
\endproclaim

\proclaim{4.8 Corollary} Let $\sE_X$ be the sheaf in (4.7). Then for an
effective divisor
$A$ in $S$ we get
$
\mu(A)^3=3(A.F)(A.A)+6\phi_1(A.F)(A.k)^2
$ with
$$\phi_1=
\delta_pq(\delta_pq+q-d)+\delta_qp(\delta_qp+p-d)+{d^2\over4}
+{d\over2}(pq-p-q).
$$
\endproclaim
\demo{Proof}
The divisor $\mu(A)$ on $X$ is the pullback along the map
$X@>>>Hilb^2(S)$ of the divisor
$$
\align
\fG=&\fA+\delta_p(\delta_p+1)(A.F_p)\fF_p+\delta_q(\delta_q+1)(A.F_q)\fF_q\\
&+{1\over4}(A.D)\left(\fD+2\delta_p\fF_p+2\delta_q\fF_q-\fK_S+\fT\right)\\
\phantom{\fG}=&\fA+x\fF+y\fT.
\endalign
$$
The fraktur--letters $\fF$ and $\fA$ denote the divisors
in $\Pic\Hilb^2(S)\cong\Pic(S)\oplus{1\over2}\bZ\fT$ corresponding to
divisors $F$ and $A$ in $S$. The divisor $\fT$ is  the exceptional locus
of the Chow map $Hilb^2(S)@>>>S^2(S)$.
So we get $\mu(A)^3=\fG^3.\fF=3(A.A)(A.F)+6x(A.F)^2-24y^2(A.F)$. The latter
formula is obtained by
applying the following multiplication rules for divisors in $\Hilb^2(S)$: \qed
\enddemo

\proclaim{4.9 Lemma} Let $S$ be an algebraic surface and $A$, $B$, $C$, $D$
divisors on $S$. Then on $Hilb^2(S)$ one has the following intersection
products:
\roster
\item $\fA.\fB.\fC.\fD=(A.B)(C.D)+(A.C)(B.D)+(A.D)(B.C)$
\item $\fA.\fB.\fC.\fT=0$
\item $\fA.\fB.\fT^2=-8(A.B)$
\item $\fA.\fT^3=8(A.K_S)$
\item $\fT^4=-8(K_S^2+c_2(S)).$
\endroster
\endproclaim
\demo{Proof} Without loss of generality, we may assume the divisors in
$S$ very ample and in general position. The first two formulas then
are obtained via the Schubfach prinziple. The last three follow from
the explicit description of the twofold cover of $\Hilb^2(S)$ as the
blow up of $S\times S$  along the diagonal.\qed
\enddemo
\medskip
\subheading{IV.c Infinitesimal extensions of bundles and sections}
\smallskip
The subspace $\sM_1^H(D,2)\subset\sM^H(D,2)$ is not reduced in
general. Let $\sM(\alpha,\beta)$ be an irreducible component
parametrizing type--1 bundles of the form
$$
0@>>>\sO_S(C)@>>>\sE@>>>\sJ_Z(D-C)@>>>0
$$
with $2C-D=-F+\alpha F_p+\beta F_q$.

\proclaim{4.10 Proposition} The generic point of $\sM(\alpha,\beta)$ has length
1 if $\alpha=\beta=0$, length 4 if $\alpha\beta\not=0$ and length 2 else.
\endproclaim

\demo{Proof} Again the proof goes along the lines of \cite{Ba2, ch. 5}, hence
we will be brief. The bundle is obtained from the unique nonsplit extension
$$
0@>>>\sO_S(C+2F)@>>>\overline\sE@>>>\sO_S(D-C)@>>>0
$$
by elementary transformation along a sheaf $\sL_I$ supported on two
disjoint generic fibers $F$ and $F'$ which restricts to line bundles of
degree one on either component. Note that the restriction of $\overline\sE$
to a generic fiber is nonsplit! (For simplicity we henceforth assume
$D-C=0$. This can be achieved by tensoring with a line bundle.)
The extension group $\Ext^1(\overline\sE,
\overline\sE)$ is zero if $\alpha$ and $\beta$ both are, is 2--dimensional,
if neither of $\alpha$ nor $\beta$ vanishes and one dimensional else. This
follows from Riemann--Roch and the computation of $\dim\Hom(\overline\sE,
\overline\sE)=3$ and $\dim\Hom(\overline\sE,\overline\sE(K_S))\in \{1,2,3\}$.
The inclusion $j:\sE\lra\overline\sE$ induces maps
$j^*:\Ext^1(\overline\sE,\overline\sE)\lra\Ext^1(\sE,\overline\sE)$ and
$j_*:\Ext^1(\sE,\sE)\lra\Ext^1(\sE,\overline\sE)$. Consider the commuting
diagram:
$$
\CD
.@>>>.@>>>\Ext^1(\overline\sE,\overline\sE)   @>>>
\Ext^1(\overline\sE,\sL_I)=0  \\
@VVV   @VVV             @V{j^*}VV             @VVV                    \\
\Hom(\sE,\sL_I)@>T>>\Ext^1(\sE,\sE)@>j_*>>\Ext^1(\sE,\overline\sE)@>>>\Ext^1(\sE,\sL_I).
\endCD
$$
Now $j^*$ is injective and the image of $T$ is the 4--dimensional
tangent space to
the reduction $\sM(\alpha,\beta)_{red}$. From (4.2) we know the dimension
of $\Ext^1(\sE,\sE)$ to be 5 (iff $\alpha\beta=0$) or 6. In particular,
$\im j^*=\im j_*$ if and only if $\alpha+\beta\not=0$. It remains to
consider two cases:

\demo{Suppose $\alpha+\beta\not=0$} Then every nonsplit extension in
$\Ext^1(\overline\sE,\overline\sE)$ lifts to a nontrivial deformation of
$\sE$ over the double point $Z_2=\Spec\bC[t]/(t^2)$ which are not tangent
to $\sM(\alpha,\beta)_{red}$. We can write down explicit representatives for
$\Ext^1(\overline\sE,\overline\sE)$
the following way:  For $\alpha\not=0$ let $\sF$ denote the nonsplit extension
$$
0\ra\sO(F+\beta F_q)\lra\sF\lra\sO\ra 0.$$
On the bundle $\sF\oplus \sF(\alpha F_p)$ define an endomorphism $\varphi$
with $\varphi^2=0$ by
$$\varphi=\left(\matrix
\tau&0\\
id\otimes\sigma&-\tau'
\endmatrix\right).
$$
Here $\sigma$ is a nontrivial section of $\sO(\alpha F_p)$ and
$\tau:\sF\ra\sO(2F+\beta F_q)\subset\sF$
is a nontrivial map chosen in such a way that the support of its
 cokernel is disjoint from $F_p$. Finally, $\tau'\in\sE nd\ \sF(\alpha F_p)$
satisfies
$\tau'\circ(id\otimes\sigma)=(id\otimes\sigma)\circ\tau$.
It is easy to check the isomorphisms
$\coker\varphi\cong\ker\varphi\cong\overline\sE$. In particular, the sheaf
$\sF\oplus\sF(\alpha F_p)$, together with the endomorphism $\varphi$ defines
an infinitesimal deformation $\overline\sE_{Z_2}$ of $\overline\sE$ over the
double point.
The infinitesimal deformation $\sE_{Z_2}$ then is the kernel of a map
of $\overline\sE_{Z_2}$ to an infinitesimal deformation of $\sL_I$.
Now suppose $\sE_{Z_3}$ extends $\sE_{Z_2}$ over the triple point
$\Spec\bC[t]/(t^3)$. Then as in the proof of \cite{Ba2, 5.1} it follows
$h^0(\sE_{Z_3}|_{pF_p})=2$. The induced map $Z_3\to\sM^H(D,2)$ thus cannot
factor through $\sM_1^H(D,2)$.
\enddemo
\demo{Suppose $\alpha+\beta=0$} Then consider the bundle $\sA_n$ over $S$
which inductively is obtained as the unique nonsplit extension
$$0@>>>\sO_S(nF)@>>>\sA_{n+1}@>>>\sA_n@>>>0$$
with $\sA_1=\sO_S$. The restriction of $\sA_n$ to each fibre is
indecomposable. Let $t: \sA_4\to\sA_4$ be a map of rank 2 factoring
through multiplication $\sA_2\to\sA_2(2F)$ with a section of $\sO_S(2F)$
vanishing at two generic fibres $F$ and $F'$. Let $x$ and $y$ be reduced points
on $F$ and $F'$. Then there is a commuting diagram of surjections
$$
\CD
\sA_4@>\varphi>>\sG=(\sH^2_F\vdash\sO_F(x))\oplus(\sH_{F'}^2\vdash\sO_{F'}(y))\\
@VVV   @VV{p_1\oplus p_2}V\\
\sA_2@>{\overline\varphi}>>\sO_F(x)\oplus\sO_{F'}(y)
\endCD
$$
of $\sO_S[t]/(t^2)$--modules (with $t:\sG\to\sG$ the zero map). Here $\sH_F^2$
denotes the unique rank--2 indecomposable bundle on $F$ with trivial
determinant
and $(A\vdash B)$ denotes a (unique) indecomposable bundle,
together with a short exact sequence
$$0@>>>A@>>>(A\vdash B)@>>>B@>>>0.$$
The kernel restricts
to $\ker\varphi|_{F\cup F'}=\sO_F(-x)\oplus\sO_{F'}(-y)\oplus \sG$.
The point now is that the map $\overline\varphi: \sH_{F\cup
F'}^2\to\sO_F(x)\oplus\sO_{F'}(y)$ extends to a map
$\sG\to\sO_F(x)\oplus\sO_{F'}(y)$.
Then let $\sE_{Z_2}$ be the kernel of the induced map
$\ker\varphi\to\sO_F(x)\oplus\sO_{F'}(y)$. It is easy to verify that
this sheaf defines a map $Z_2\to\sM^H(D,2)$ which does not factor through
the subscheme of 1--instantons.
\qed\enddemo\enddemo

Using this explicit description of infinitesimal deformations of bundles
it is easy to find which sections extend to these deformations:

\proclaim{4.12 Corollary} Let $\xi\in\sM(\alpha,\beta)$ denote the generic
point and $\sE_\xi$ the pull back of a universal bundle to $S\times \xi$.
 Then
$\length \sE xt^2_{p_{\xi}}(\sE_\xi,p_S^*\omega_S)$ is
given by the sum $\sum_{X\in\Delta}h^0(S;\sO_S(X))$
over the set $\Delta=
\{C,(C-\alpha F_p), (C-\beta F_q), (C-\alpha F_p-\beta F_q)\}$.
\qed\endproclaim

Now let $\sN(\alpha,\beta)$ denote an irreducible component in $\sM^1(D,2)$
parametrizing type--2 bundles generically of the form
$$
0@>>>\sO_S(C)\oplus\sO_S(C-F+\alpha F_p+\beta F_p)@>>>\sE@>>>\sJ_{Z\subset
F}@>>>0
$$
with $0<\alpha,\beta$ and $\alpha q+\beta p<pq$. One easily computes
$\Ext^2(\sE,\sE)=0$ and therefore $\sM^H(D,2)$ is generically smooth
along $\sN(\alpha,\beta)$.

\proclaim{4.13 Proposition} The generic point of $\sN(\alpha,\beta)$
has length 1.\endproclaim

\demo{Proof} It suffices to construct an infinitesimal deformation $\sE_{Z_2}$
such that no section of $\sE$ extends (compare the proof of 4.11).
Let $t:\sA_2\to\sA_2$ be a rank--1 map factoring through multiplication
 with a section of $\sO_S(F)$ vanishing at $F$ and set
$\sJ_{Z\subset F}^\vee=\sH om(\sJ_{Z\subset F},\sO_F)$. Then there is a
commuting diagram
$$
\CD
\sA_2\oplus\sA_2@>\varphi>>(\sH^2_F\oplus\sH^2_F)\vdash\sJ_{Z\subset F}^\vee\\
@VVV   @VVV\\
\sA_1\oplus\sA_1@>{\overline\varphi}>>\sJ_{Z\subset F}^\vee
\endCD
$$

of $\sO_S[t]/t^2$--modules. The kernel restricts to $\ker\varphi|_F=
\sJ_{Z\subset F}\oplus\sJ_{Z\subset F}^\vee$ and the map
$\overline\varphi:\sH_F^2\oplus\sH_F^2\to\sJ_{Z\subset F}^vee$ extends to a
map $(\sH^2_F\oplus\sH^2_F)\vdash\sJ_{Z\subset F}^\vee\to\sJ_{Z\subset
F}^\vee$.
The sheaf $\sE_{Z_2}$ then is the kernel of the induced map $\ker\varphi\to
\sJ_{Z\subset F}^\vee$.\qed
\enddemo

\proclaim{4.14 Corollary} If $\xi\in\sN(\alpha,\beta)$ denotes the generic
point and $\sE_\xi$ the pull back of a universal bundle to $S\times \xi$.
 Then
$\length \sE xt^2_{p_{\xi}}(\sE_\xi,p_S^*\omega_S)$ is
$h^0(S;\sO_S(C)\oplus\sO_S(C_F+\alpha F_p+\beta F_p))$.\qed
\endproclaim
\bigskip
 \bigskip
\heading
V. Combinatorics
\endheading
We finally combine the information of the various strata to compute the
first two coefficients in the $\Spin^\bC $--invariant for Dolgachev surfaces.
The result is:

\proclaim{5.1 Theorem} The first two coefficients in the
polynomial $q_S(n)$ are  $a(n)=3n$ and $b(n)=(2p^2q^2-2p^2-2q^2-1)n$.
\endproclaim

\demo{Remarks} - Here we use the symmetrization convention which gives
$$
(Q^2+Qk^2+k^4)(A,A,A,A)=(A.A)^2+(A.A)(A.k)^2+(A.k)^4,
$$
which actually differs from the one used in \cite{Ba2}.\hfill\break
- In the case $p=q=1$ one can actually determine the third coefficient to
be $c(n)=21n$.
 Our methods, however, do not allow to compute the third coefficient in
general. Such a computation probably
would give the diffeomorphism classification of Dolgachev surfaces,
which Friedman (\cite{Fr})
 recently obtained by other means. Regarding the problem of diffeomorphism
classification, the above coefficients contain exactly the same information
as the $\Gamma$--invariant of Donaldson (cf. \cite{Ba1}).
\enddemo

The above convention gives $Q^2(A,A,A,F)=(A.A)(A.F)$ and
$Qk^2(A,A,A,F)={1\over2}(A.F)(A.k)^2$. To  prove the theorem, we
first gather the contributions of the different strata into a closed formula.
We will assume $p$ to be odd. The map
$$
(\alpha,\beta)\mapsto\cases
\left(\tfrac{\alpha-1}2,
\tfrac{\beta-1}2\right)&\text{if $\alpha$ is odd}\\
\left(\tfrac{p-\alpha-1}2,\tfrac{q-\beta-1}2\right)&\text{if $\alpha$ is even}
\endcases
$$
defines a bijection of
$$\square= \{(\sigma,\tau)\in\bZ^2\,\vert\,0\le\sigma\le\tfrac{p-1}2,\,
0\le\tau\le\tfrac{q-1}2\}$$
with the set
$$\{(\alpha,\beta)\,\vert\,(\alpha-1)q+(\beta-1)p\,\,\,\text{even},\,
\tfrac\alpha p+\tfrac\beta q<1\}$$
indexing the type--1
components $\sM(\alpha,\beta)$. Similarly, the map
$$
(\alpha,\beta)\mapsto\cases
\left(\tfrac{\alpha-1}2,
\tfrac{q-\beta-1}2\right)&\text{if $\alpha$ is odd}\\
\left(\tfrac{p-\alpha-1
}2,\tfrac{\beta-1}2\right)&\text{if $\alpha$ is even}
\endcases
$$
defines a bijection of $\square$ with the indexing set
$$\{(\alpha,\beta)\,\vert\,\alpha q+\beta p\,\,\,\text{odd},\,
\tfrac\alpha p+\tfrac\beta q<1\}$$
of type--2 components $\sN(\alpha,\beta)$.
(Note that for $\alpha=0$ or $\beta=0$ there is no component
$\sN(\alpha,\beta)$. But then the formula (4.5) also vanishes.)

Next come multiplicities. These multiplicities merely count the number of
sections or cosections of a generic bundle parametrized in $\sM(\alpha,\beta)$
or $\sN(\alpha,\beta)$. We fix a number $n=lpq+Aq+Bp$ with our usual
convention $0\le
{A\over p},{B\over q}<1$.

\proclaim{5.2 Lemma}
The multiplicity of
the component $\sM(\alpha,\beta)$
or
 $\sN(\alpha,\beta)$ is given by
$$m(\sigma,\tau,n)=
\left(l+1+H_p(\sigma,A)+
H_q(\tau,B)\right)c(\sigma,\tfrac{p-1}2)
c(\tau,\tfrac{q-1}2)$$
$$\text{with}\,\,\,\, H_p(\sigma,A)=\cases
{1\over2}&\text{if $\sigma\ge p-A$}\\
-{1\over2}&\text{if $\sigma\ge A$}\\
0&\text{else}
\endcases
$$
$$\text{and}\,\,\,\,c(\sigma,\tfrac{p-1}2)=\cases
1&\text{if $\sigma=\tfrac{p-1}2$}\\
2&\text{else.}
\endcases
$$
\endproclaim

\demo{Proof} Let $\sE$ be a generic bundle in $\sM(\alpha,\beta)$. Then
because of (3.14)  one has to compute $$\sum_Ih^0(\sE(X))-h^2(\sE(X))$$
with $I=\{0,-\alpha F_p, -\beta F_q, -\alpha F_p -\beta F_q\}$.
One verifies the claim via checking the different cases: whether or not
$2A\le p$, whether or not $2B\le q$ and whether or not $\alpha $ is odd.
For example, if $2A\le p$, $2B\le q$ and  $\alpha $ odd, one has a presentation
$$
0\ra\sO(lF+(A+\sigma)F_p+(B+\tau)F_q)@>>>\sE@>>>
\sJ_{Z}(lF+(A-\sigma)F_p+(B-\tau)F_q+K_S)\ra 0
$$
and the claim follows by applying (3.1).

For the components $\sN(\alpha,\beta)$ the divisor $\mu(F)$ already has
multiplicity 2 (compare IV.b) which has to be multiplied with
$h^0(\sE)-h^2(\sE)$.
Again one applies (3.1) to $\sE$, which in the example above has a presentation
$$
\multline
0@>>>
\sO_S(lF+(A+\sigma)F_p+(B-\tau-1)F_q)\oplus\sO_S(lF+(A-\sigma-1)F_p+(B+\tau)F_q)
\\@>>>\sE@>>>\sJ_{Z\subset F}@>>>0.\qed
\endmultline
$$
\enddemo

The correction term
$s_q=s_q(\tau)=d(\sigma,\tau)-\left[\dfrac{d(\sigma,\tau)}{2q}\right]2q$
takes odd values in the interval $0\le s_q(\tau)<2q$.
Using this correction term we get the following
expression for the function $\phi_1$:
$$
\phi_1=\frac14(2dpq-d^2)-\frac14(2qs_q-s_q^2)-
\frac14(2ps_p-s_p^2)
$$
where $d$ as a function of $\sigma$ and $\tau$ is given by
$$
d=\cases
pq-(2\sigma+1)q-(2\tau+1)p&\text{if $2\sigma q+2\tau p<pq-p-q$}\\
-pq+(2\sigma+1)q-(2\tau+1)p&\text{else.}
\endcases
$$
The contribution of type--2 bundles is given by
$$
\phi_2=\cases
(2\sigma+1)(2\tau+1)pq&\text{if $2\sigma q+2\tau p<pq-p-q$}\\
(p-2\sigma-1)(q-2\tau-1)pq&\text{else.}
\endcases
$$
To get decent formulae, define the functions
$$
\align
T_p=&T_p(\sigma)=(2\sigma+1)(2p-2\sigma-1)\\
S_p=&S_p(\sigma)=2ps_p-s_p^2\\
R=&R_p(\sigma,\tau)=\max(pq-p-q-2\sigma q-2\tau p,0)
\endalign
$$
and their analogs $T_q$ and $S_q$ by interchanging the r\^oles
of $(p,\sigma)$ and $(q,\tau)$. Then we can sum up everything to get the

\proclaim{5.3 Lemma} The coefficients $a(n) $ and $b(n)$ are
given by the formulae
$$
\align
a(n)&=\sum_\square m(\sigma,\tau,n)\\
b(n)&=\sum_\square m(\sigma,\tau,n)\Phi(\sigma,\tau)\\
\endalign
$$
with
$\Phi=\Phi(\sigma,\tau)=12\phi_1+6\phi_2=3(q^2T_p-S_p)-3(p^2T_q+S_q)+
6pqR-3p^2q^2$.
\endproclaim

The theorem now follows from lemmata (5.4)--(5.7) below:

\proclaim{5.4 Lemma}
$$\sum_\square m(\sigma,\tau,n)=n$$
\endproclaim
\demo{The proof} is left to the reader.\qed\enddemo

\proclaim{5.5 Lemma}
$$3\sum_\square mS_q=n(2q^2+1)+
\cases
-Bp(2q^2+1)+3p\underset {0\le\tau<B}\to\sum S_q&\text{if $2B\le q$}\\
(q-B)p(2q^2+1)+3p\underset{0\le\tau<q-B}\to\sum S_q&\text{else.}
\endcases
$$
\endproclaim

\proclaim{5.6 Lemma}
$$3\sum_\square m(p^2T_q-S_q)=n(2q^2+1)(p^2-1)+X_q$$
$$\text{with}\,\,\,\,\,\,\,
X_q=\cases
-Bp(2q^2+1)(p^2-1)+3p\underset {0\le\tau<B}\to\sum (p^2T_q-S_q)&\text{if $2B\le
q$}\\
(q-B)p(2q^2+1)(p^2-1)+3p\underset {0\le\tau<q-B}\to\sum
(p^2T_q-S_q)&\text{else.}
\endcases
$$
\endproclaim

\demo{Proof of 5.5}
For $0\le\tau<{q\over2}$ each term $S_q(\tau)$ is a product $x(2q-x)$
with $0<x\le q$ an odd integer and each such value is obtained exactly
once. So we get
$$
6\sum_{0\le\tau<\left[\tfrac{q}2\right]}S_q=
\cases
q(q-1)(2q-1)&\text{if $q$ is odd}\\
q(2q^2+1)&\text{else.}
\endcases
$$
In particular the "$(l+1)cc$"--summand adds up to
$$
3\sum_\square(l+1)c(\sigma,\tfrac{p-1}2)c(\tau,\tfrac{q-1}2)S_q=(l+1)pq(2q^2+1).$$
Similarly one gets the two other components
$$
\align
6\sum_\square H_p(\sigma,A)c(\sigma,\tfrac{p-1}2)c(\tau,\tfrac{q-1}2)S_q&=
(2A-p)q(2q^2+1)\\
6\sum_\square H_q(\tau,B)c(\sigma,\tfrac{p-1}2)c(\tau,\tfrac{q-1}2)S_q&=
\cases
-pq(2q^2+1)+6p\underset {0\le\tau<B}\to\sum S_q&\text{if $2B\le q$}\\
pq(2q^2+1)-6p\underset {0\le\tau<q-B}\to\sum S_q&\text{else.\,\qed}
\endcases
\endalign
$$
\enddemo

\demo{Proof of 5.6} The values $S_q(\tau)$ are merely permutations
of the values of $T_q(\tau)$ for $0\le\tau<[\tfrac{q}2]$ and
$S_q(\tfrac{q-1}2)=T_q(\tfrac{q-1}2)$. So the claim is immediate from
(3.5).\qed
\enddemo

\proclaim{5.7 Lemma}
$$6pq\sum_\square mR=n(p^2-1)(q^2-1)-X_p-X_q.$$
\endproclaim

\demo{Proof} For $2\sigma q+2\tau p<pq-p-q$ the product
$c(\sigma,\tfrac{p-1}2)c(\tau,\tfrac{q-1}2)=4$ is constant on the index set.
So we get
$$
\sum_\square mR=\sum_\square 4(l+1)R +\Xi_A+\Xi_B
$$
$$\text{with}\,\,\,\,\,\,\,\Xi_B=\cases
\sum_\square(-2R)+\underset {0\le\tau<B}\to\sum \underset{0\le\sigma<\tfrac
q2}\to\sum 2R&\text{if $2B\le q$}\\
\sum_\square2R-\underset {0\le\tau<q-B}\to\sum \underset{0\le\sigma<\tfrac
q2}\to\sum (-2R)&\text{else}
\endcases
$$
and an analog expression for $\Xi_A$. Setting
$z(\tau)=\left[\tfrac{pq-p-q-2\tau p}{2q}\right]$ we get
$$
\align
\sum_{0\le\tau<y}\sum_{0\le\tau<\tfrac q2}R
&=
\sum_{0\le\tau<y}\sum_{0\le\sigma<z(\tau)}(pq-p-q-2\sigma q-2\tau p)\\
&=
\underset {0\le\tau<B}\to\sum \left(z(\tau)+1\right)
\left(pq-p-q-2\tau p -qz(\tau)\right)\\
&=\underset {0\le\tau<B}\to\sum
\frac1{4p}\left((q^2(p^2-1)-(p^2T_q-S_q)\right).
\endalign
$$
Using again (5.5), we may replace $y$ by the different choices we have and get
$$6pq\sum_\square4(l+1)R=(l+1)(p^2-1)(q^2-1),$$
$$
\Xi_B=\cases
(-\tfrac1{12}(q^2-1)+\tfrac q4B)(p^2-1)-\underset {0\le\tau<B}\to\sum
\tfrac1{4q}(p^2T_q-S_q)&\text{if $2B\le q$}\\
(\tfrac1{12}(q^2-1)-\tfrac q4(q-B))(p^2-1)+\underset {0\le\tau<q-B}\to\sum
\tfrac1{4q}(p^2T_q-S_q)&\text{else.}
\endcases
$$
The result follows.\qed
\enddemo

\enddocument